\documentclass[12pt]{article}
\usepackage[utf8]{inputenc}
\usepackage{tikz}
\usepackage{amssymb,amsmath,caption,subcaption,graphicx,hyperref, xcolor}
\usepackage{amsthm}
\usepackage{authblk, blindtext}
\usepackage{lineno}
\usetikzlibrary{decorations.pathreplacing}
\usepackage{bbm}

\usepackage{enumitem}

\usepackage{lipsum}

\newenvironment{stretchpars}
 {\par\setlength{\parfillskip}{0pt}}
 {\par}

\usepackage{tkz-graph}
\tikzstyle{vertex}=[circle, draw, inner sep=3pt, minimum size=6pt]

\usepackage{tikzsymbols}
\usepackage{textcomp}
\usepackage{parskip}

\usepackage{pgfplots}
\usepackage{verbatim}
\usepackage{lineno} 
\usepackage{subfiles}
\usepackage{xcolor}
\usepackage{blindtext}

\newtheorem{thm}{Theorem}[section]

\newtheorem{prop}[thm]{Proposition}
\newtheorem{cor}[thm]{Corollary}
\newtheorem{lem}[thm]{Lemma}

\newtheorem{conj}[thm]{Conjecture}
\newtheorem{ex}[thm]{Example}

\newtheorem{obs}[thm]{Observation}

\definecolor{purple}{rgb}{0.64 ,0.17, 1.0}
\definecolor{MyGreen}{rgb}{0.15, 0.5, 0.06}
\definecolor{BurntOrange}{rgb}{0.7,0.35,0.0}
\definecolor{pink}{rgb}{1.0, 0.1, 0.5}

 \newcommand{\bpf}{\begin{proof}}
 \newcommand{\epf}{\end{proof}}

 \newcommand{\C}{\mathcal{C}}
 \newcommand{\aw}{\textup{aw}}
 \newcommand{\dist}{\textup{d}}
 \newcommand{\diam}{\textup{diam}}

\definecolor{cbblue}{RGB}{68,119,170}
\definecolor{cbcyan}{RGB}{102,204,238}
\definecolor{cbgreen}{RGB}{34,136,51}
\definecolor{cbyellow}{RGB}{204,187,68}
\definecolor{cbred}{RGB}{238,102,119}
\definecolor{cbpurple}{RGB}{170,51,119}
\definecolor{cbgray}{RGB}{187,187,187}

\title{Anti-van der Waerden Numbers of Graph Products with Trees}

\author[1]{Zhanar Berikkyzy}
\author[2]{Joe Miller}
\author[3]{Elizabeth Sprangel}
\author[4]{Shanise Walker}
\author[5]{Nathan Warnberg}

\affil[1]{Mathematics Department, Fairfield University\{zberikkyzy@fairfield.edu\}}
\affil[2]{Department of Mathematics, Iowa State University \{jmiller0@iastate.edu\}}
\affil[3]{Health and Human Services Department, Hennepin County, \{elizabeth.sprangel@hennepin.us\}}
\affil[4]{Department of Mathematical Sciences, Clark Atlanta University, \{swalker@cau.edu\}}
\affil[5]{Department of Mathematics and Statistics, University of Wisconsin-La Crosse, \{nwarnberg@uwlax.edu\}}

\begin{document}

\maketitle
\begin{abstract}
    Given a graph $G$, an exact $r$-coloring of $G$ is a surjective function $c:V(G) \to [1,\dots,r]$.  An arithmetic progression in $G$ of length $j$ with common difference $d$ is a set of vertices $\{v_1,\dots, v_j\}$ such that $\dist(v_i,v_{i+1}) = d$ for $1\le i < j$.  An arithmetic progression is rainbow if all of the vertices are colored distinctly.  The fewest number of colors that guarantees a rainbow arithmetic progression of length three is called the anti-van der Waerden number of $G$ and is denoted $\aw(G,3)$.  It is known that $3 \le \aw(G\square H,3) \le 4$.  Here we determine exact values $\aw(T\square T',3)$ for some trees $T$ and $T'$, determine $\aw(G\square T,3)$ for some trees $T$, and determine $\aw(G\square H,3)$ for some graphs $G$ and $H$.
\end{abstract}

\section{Introduction}

Ramsey Theory is the study of choosing a set of mathematical elements $S$, then coloring each element via a surjective coloring function $c:S \to \{1,2,\dots,r\}$ and determining if there are substructures that are monochromatic.  
Despite the name, Ramsey Theory was first studied by Schur in 1917 when he showed that we can always find an $N$ such that if $n \ge N$ and we color $[n] = \{1,2,\dots,n\}$ with $r$ colors, there must always be a monochromatic solution to $x+y=z$ (see \cite{S}).  Van der Waerden, in 1927, proved that given a fixed $r$ and $k$, there exists an $N$ such that if $n\ge N$ and you color $[n]$ with $r$-colors, there must a monochromatic $k$-term arithmetic progression (see \cite{W27}).  
Finally, Ramsey, in 1928, showed that given a fixed $r$ and $k$, there exists an $N$ so that if $n \ge N$, then no matter how you color the edges of the complete graph $K_n$ using $r$ many colors, there must be a monochromatic complete subgraph $K_k$ (see \cite{R}).  
It was not until 1973, in \cite{ESS}, when Erd{\H{o}}s, Simonovits, and S{\'{o}}s introduced the idea of looking for polychromatic or rainbow structures.  
Looking for rainbow structures are known as Anti-Ramsey problems, or in the case of this paper, Anti-van der Waerden problems.

A \emph{$k$-term arithmetic progression}, $k$-AP, in $[n]$ is of the form
\[ a, a+d, a+2d, \dots, a+(k-1)d\]
where $k \ge 2$ and $d\ge 1$.  
An \emph{exact $r$-coloring} of $[n]$ is a surjective function $c:[n] \to [r]$.

The anti-van der Waerden number of $[n]$ with respect to $k$-APs, first defined in~\cite{U} in $2013$ and denoted $\aw([n],k)$, is the smallest $r$ such that every exact $r$-coloring of $[n]$ guarantees a rainbow arithmetic $k$-AP.  
Note that rainbow $3$-APs had been studied previously when balanced colorings are used in~\cite{AF,AM,J}.  
It should also be noted that $\aw([n],3)$ was studied in~\cite{DMS} and the exact value was determined in~\cite{BSY}.  
Additionally, rainbow solutions to linear equations in $[n]$ were studied in~\cite{budden} and~\cite{FGRWW}.
Interestingly, the authors in~\cite{budden} and~\cite{FGRWW} produced similar results simultaneously but using differing arguments.
The definition of $\aw([n],k)$ was also extended to $\mathbb{Z}_n$ to get $\aw(\mathbb{Z}_n,k)$ in~\cite{DMS} which led to results on arithmetic progression in finite abelian groups by Young in~\cite{finabgroup}.  
Rainbow solutions to linear equations in $\mathbb{Z}_n$ have also been studied in~\cite{BKKTTY,RFC,LM}.  
Since $[n]$ behaves like a path graph and $\mathbb{Z}_n$ behaves like a cycle, anti-van der Waerden numbers in graphs, particularly tree graphs, and graph products were studied in~\cite{SWY, MW, RSW}.  
This paper continues in this vein of research.

Graphs $G=(V,E)$ in this paper are simple, undirected and connected so edge $\{u,v\}$ will be shortened to $uv\in E(G)$ where $u,v$ are vertices in the vertex set $V(G)$.   If $uv\in E(G)$, then $u$ and $v$ are \emph{neighbors} of each other. 
The \emph{distance} between vertex $u$ and $v$ in graph $G$ is denoted $\dist_G(u, v)$, or just $\dist(u, v)$ when context is clear, and is the smallest length of any $u-v$ path in $G$.
A $u-v$ path of length $\dist(u,v)$ is called a $u-v$ \emph{geodesic}.
Graph $G'$ is a \emph{subgraph} of $G$ if $V(G') \subseteq V(G)$ and $E(G') \subseteq E(G)$. A subgraph $G'$ of $G$ is an \emph{induced subgraph} if whenever $u$ and $v$ are vertices of $G'$ and $uv$ is an edge of $G$, then $uv$ is an edge of $G'$.  If $S$ is a nonempty set of vertices of $G$, then the \emph{subgraph of $G$ induced by $S$} is the induced subgraph with vertex set $S$ and is denoted $G[S]$.  An \emph{isometric subgraph} $G'$ of $G$ is a subgraph such that $\dist_{G'}(u,v) = \dist_G(u,v)$ for all $u,v \in V(G')$.  If $G = (V,E)$ and $H = (V', E')$ then the \emph{Cartesian
product}, written $G\square H$, has vertex set $\{(x, y) : x \in V \text{ and } y \in V' \}$ and $(x, y)$ and
$(x', y')$ are adjacent in $G \square H$ if either $x = x'$ and $yy' \in E'$ or $y = y'$ and $xx' \in E$. This paper will use the convention that if \[V(G) = \{u_1,\dots, u_{n_1}\} \quad \text{and} \quad V(H) = \{w_1,\dots,w_{n_2}\},\] then $V(G\square H) = \{v_{1,1},\dots, v_{n_1,n_2}\}$ where $v_{i,j}$ corresponds to the vertices $u_i \in V(G)$ and $w_j \in V(H)$.  Also, if $1\leq i \leq n_2$, then $G_i$ denotes the $i$th labeled copy of $G$ in $G \square H$. Likewise, if $1 \leq j \leq n_1$, then $H_j$ denotes the $j$th labeled copy of $H$ in $G \square H$.  In other words, $G_i$ is the induced subgraph $G_i = G\square H[\{v_{1,i},\dots, v_{n_2,i}\}]$, and $H_j$ is the induced subgraph $H_j = G\square H[\{v_{j,1}, \dots, v_{j,n_1}\}]$.  Notice that the $i$ subscript in $G_i$ corresponds to the $i$th vertex of $H$ and the $j$ in the subscript in $H_j$ corresponds to the $j$th vertex of $G$ ( see Example \ref{ex:cartprod} below).

\begin{stretchpars}A \emph{$k$-term arithmetic progression in graph $G$} ($k$-AP) is a set of vertices \end{stretchpars}\vspace{-.1in}$\{v_1,\dots,v_k\}$ such that $\dist(v_i,v_{i+1}) = d$ for all $1\le i \le k-1$.  
A $k$-term arithmetic progression is \emph{degenerate} if $v_i = v_j$ for any $i\neq j$.  
Note that technically, since a $k$-AP is a set, the order of the elements does not matter.  
However, oftentimes $k$-APs will be presented in the order that provides the most intuition. 
An \emph{exact $r$-coloring of a graph $G$} is a surjective function $c:V(G) \to [r]$.  
A set of vertices $S$ is \emph{rainbow} under coloring $c$ if for every $v_i,v_j\in S$, $c(v_i) \neq c(v_j)$ when $i\neq j$.  
Given a set $S\subset V(G)$, define $c(S) = \{ c(s) \, | \, s\in S\}$.
The \emph{anti-van der Waerden number of graph $G$ with respect to $k$}, denoted $\aw(G,k)$, is the least positive number $r$ such that every exact $r$-coloring of $G$ contains a rainbow $k$-term arithmetic progression.  
If $|V(G)| = n$ and no coloring of the vertices yields a rainbow $k$-AP, then $\aw(G,k) = n+1$. 

In this paper, $P_n$ denotes that path graph on $n$ vertices and $C_n$ denotes the cycle graph on $n$ vertices. 
Example~\ref{ex:cartprod} below discusses the anti-van der Waerden number of $P_3\square C_6$.

\begin{ex}\label{ex:cartprod}

    Consider the graph $P_3\square C_6$ where $V(P_3) = \{u_1,u_2,u_3\}$ and $V(C_6) =\{w_1,w_2,w_3,w_4,w_5,w_6\}$.  Let $G=P_3$ and $H = C_6$ as in the definition.  Now, $G_5$ is a subgraph of $P_3\square C_6$ that is isomorphic to $P_3$ and corresponds to vertex $w_5$ of $C_6$.  Similarly, $H_2$ is a subgraph of $P_3\square C_6$ that is isomorphic to $C_6$ and corresponds to vertex $u_2$ of $P_3$.
    Note that coloring $c(v_{1,1}) = red$, $c(v_{3,4}) = blue$ and all other vertices $white$ avoids rainbow $3$-APs (see Figure \ref{fig:cartprodex} below).  Since we have a rainbow-free, exact $3$-coloring this implies that $4\le\aw(P_3\square C_6,3)$.  Some case analysis, or the use of Theorem \ref{theorem:rsw}, gives that $\aw(P_3\square C_6,3)\le 4$. Therefore, $\aw(P_3\square C_6,3) = 4$.
  
  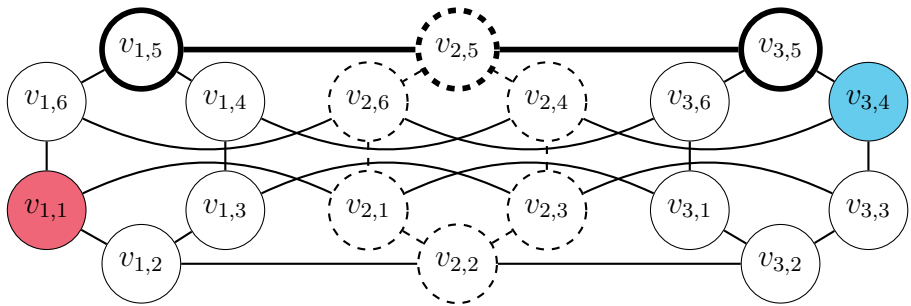
\begin{figure}[ht!]
      \centering

    \begin{tikzpicture}[scale = .95]
	    \node[draw,circle, fill = cbred] (11) at (0,.75) {$v_{1,1}$};
		\node[draw,circle] (12) at (1.325,0) {$v_{1,2}$};
		\node[draw,circle] (13) at (2.5,.75) {$v_{1,3}$};
		\node[draw,circle] (14) at (2.5,2.27) {$v_{1,4}$};
		\node[draw,circle, line width=0.75mm] (15) at (1.325,3) {$v_{1,5}$};
		\node[draw,circle] (16) at (0,2.27) {$v_{1,6}$};
		\node[draw,circle,dashed, thick] (21) at (4.5,.75) {$v_{2,1}$};
		\node[draw,circle,dashed, thick] (22) at (5.75,0) {$v_{2,2}$};
		\node[draw,circle,dashed, thick] (23) at (7,.75) {$v_{2,3}$};
		\node[draw,circle,dashed, thick] (24) at (7,2.27) {$v_{2,4}$};
		\node[draw,circle,dashed, line width=0.75mm] (25) at (5.75,3) {$v_{2,5}$};
		\node[draw,circle,dashed, thick] (26) at (4.5,2.27) {$v_{2,6}$};
	    \node[draw,circle] (31) at (9,.75) {$v_{3,1}$};
		\node[draw,circle] (32) at (10.27,0) {$v_{3,2}$};
		\node[draw,circle] (33) at (11.5,.75) {$v_{3,3}$};
		\node[draw,circle, fill = cbcyan] (34) at (11.5,2.27) {$v_{3,4}$};
		\node[draw,circle, line width=0.75mm] (35) at (10.27,3) {$v_{3,5}$};
		\node[draw,circle] (36) at (9,2.27) {$v_{3,6}$};
		
      	\draw[thick]  (11) to node [auto] {} (12);
      	\draw[thick]  (12) to node [auto] {} (13);
      	\draw[thick]  (13) to node [auto] {} (14);
      	\draw[thick]  (14) to node [auto] {} (15);
      	\draw[thick]  (15) to node [auto] {} (16);
      	\draw[thick]  (16) to node [auto] {} (11);
      	\draw[thick,dashed]  (21) to node [auto] {} (22);
      	\draw[thick,dashed]  (22) to node [auto] {} (23);
      	\draw[thick,dashed]  (23) to node [auto] {} (24);
      	\draw[thick,dashed]  (24) to node [auto] {} (25);
      	\draw[thick,dashed]  (25) to node [auto] {} (26);
      	\draw[thick,dashed]  (26) to node [auto] {} (21);
      	\draw[thick]  (31) to node [auto] {} (32);
      	\draw[thick]  (32) to node [auto] {} (33);
      	\draw[thick]  (33) to node [auto] {} (34);
      	\draw[thick]  (34) to node [auto] {} (35);
      	\draw[thick]  (35) to node [auto] {} (36);
      	\draw[thick]  (36) to node [auto] {} (31);
      	\draw[thick, bend left = 25]  (11) to node [auto] {} (21);
      	\draw[thick]  (12) to node [auto] {} (22);
      	\draw[thick, bend left = 25]  (13) to node [auto] {} (23);
      	\draw[thick, bend right = 25]  (14) to node [auto] {} (24);
      	\draw[thick, line width=0.75mm]  (15) to node [auto] {} (25);
        \draw[thick, bend right = 25]  (16) to node [auto] {} (26);
      	\draw[thick, bend left = 25]  (21) to node [auto] {} (31);
      	\draw[thick]  (22) to node [auto] {} (32);
      	\draw[thick, bend left = 25]  (23) to node [auto] {} (33);
      	\draw[thick, bend right = 25]  (24) to node [auto] {} (34);
      	\draw[thick, line width=0.75mm]  (25) to node [auto] {} (35);
        \draw[thick, bend right = 25]  (26) to node [auto] {} (36);
    \end{tikzpicture}
     
      \caption{The subgraph $G_5$ has been bolded and $H_2$ has been dashed. There is also an exact $3$-coloring using the colors $red$, $white$ and $blue$ that avoids rainbow $3$-APs.} 
      
      \label{fig:cartprodex}
  \end{figure}
  
\end{ex}

\section{Preliminary Results}

Proposition \ref{proposition:dist} will be used regularly without reference as the result is rather intuitive.
\begin{prop}\label{proposition:dist}\cite{MW}
If $v_{i,j},v_{h,k} \in V(G \square H)$, then \[\dist_{G\square H}(v_{i,j},v_{h,k}) = \dist_G(u_i,u_h) + \dist_H(w_j,w_k).\]
\end{prop}

The diameter of graph $G$, denoted $\diam(G)$, is the maximum distance between all pairs of vertices in $V(G)$.

\begin{obs}\label{obs:diamGxH}
If $G$ and $H$ are graphs, then a direct consequence of Proposition \ref{proposition:dist} is that \[\diam(G\square H) = \diam(G) + \diam(H).\] 
\end{obs}

\begin{cor}\label{cor:isosubprod}\cite{MW}
    If $G'$ is an isometric subgraph of $G$ and $H'$ is an isometric subgraph of $H$, then $G'\square H'$  is an isometric subgraph of $G\square H$. 
\end{cor}

Lemma \ref{isometricpathorC3} is powerful since it guarantees isometric subgraphs.  Isometric subgraphs are important when investigating anti-van der Waerden numbers because distance preservation implies $k$-AP preservation. Lemma \ref{|c(V(Gi U Gj))|<3} helps restrict the number of colors each copy of $G$ or $H$ can have within $G\square H$. 

\begin{lem}\cite{RSW}\label{isometricpathorC3}
If $G$ is a connected graph on at least three vertices with an exact $r$-coloring
$c$ where $r \ge 3$, then there exists a subgraph $G'$ in $G$ with at least three colors
where $G'$ is either an isometric path or $G' = C_3$.
\end{lem}

Using Corollary \ref{cor:isosubprod} and Lemma \ref{isometricpathorC3} we obtain Theorem \ref{theorem:GxH3}. 

\begin{thm}\label{theorem:GxH3}
Let $G$ be a connected graph on at least three vertices. If $diam(G) \le 2$ and $aw(P_3\square H,3) = 3$, then $aw(G\square H,3)=3$.
\end{thm}
\bpf
 
Let $c$ be an exact $3$-coloring of $G\square H$. 
By Lemma \ref{isometricpathorC3}, there exists a subgraph $X$ of $G\square H$ with $3$-colors such that $X$ is either an isometric path or is isomorphic to $C_3$. 
If $X\cong C_3$, this forms a rainbow $3$-AP and we are done. 
So suppose that $X$ is an isometric path in $G\square H$. 
Notice that X is a subgraph of $P_3\square H$, which in turn is a subgraph of $G\square H$.

Notice that $P_3\square H$ is $3$-colored under $c$, because it contains an isometric path. 
Since $aw(P_3\square H,3) = 3$, the graph $P_3\square H$ contains a rainbow $3$-AP under the coloring $c$. 
Since $P_3\square H$ is an isometric subgraph of $G\square H$, by Corollary \ref{cor:isosubprod}, $G\square H$ contains this rainbow $3$-AP as well. 
Therefore, $aw(G\square H,3)=3$. 
\epf

Examples of graphs $G$ that satisfy Theorem~\ref{theorem:GxH3} are the complete graph on $n$ vertices $K_n$, the complete bipartite graph on $m$ and $n$ vertices $K_{m,n}$, or graphs with a dominating vertex. Examples of graphs $H$ that satisfy $aw(P_3\square H,3) = 3$  are $P_n$ with odd $n$, $C_{2k+1}$, $C_{4k}$ (see \cite{MW}) or so-called $3$-peripheral trees (see Section \ref{3peripheraltrees} for a discussion of $3$-peripheral trees). Thus, we can cross any graphs $G$ and $H$ above to get an anti-van der Waerden number of three.

\begin{lem}\cite{RSW}\label{|c(V(Gi U Gj))|<3}
Assume $G$ and $H$ are connected with $|V(H)| \geq 3$. Suppose $c$ is an exact, rainbow-free $r$-coloring of $G\square H$, such that $r \geq 3$ and $|c(V(G_i))| \leq 2$ for $1 \leq i \leq n$. If $w_iw_j \in E(H)$, then
$|c(V(G_i) \cup V (G_j))| \leq 2$.
\end{lem}

Lemmas \ref{|c(Hi)|<3} and \ref{c(H_i)/c(H_j) < 2} are structural results that significantly limit the way we can color $G\square H$ and be rainbow-free.

\begin{lem}\label{|c(Hi)|<3}\cite{MW}
If $G$ and $H$ are connected, $|G|,|H| \ge 2$ and $c$ is an exact $r$-coloring of $G\square H$, $3\le r$, that avoids rainbow $3$-APs, then $|c(V(G_i))| \leq 2$ for $1 \leq i \leq |H|$.
\end{lem}

\begin{lem}\label{c(H_i)/c(H_j) < 2}\cite{SWY}
Let $G$ be a connected graph on $m$ vertices and $H$ be a connected graph on $n$ vertices. Let $c$ be an exact $r$-coloring of $G\square H$ with no rainbow $3$-APs. If $G_1,G_2, \ldots,G_n$ are the labeled copies of $G$ in $G\square H$, then $|c(V(G_j)) \setminus
c(V(G_i))| \leq 1$ for all $1 \leq i, j \leq n$.
\end{lem}

Theorem \ref{theorem:rsw} was one of the results that initiated the vein of research that this paper is following.  In particular, this paper is trying to characterize when the product of two graphs has an anti-van der Waerden number of three and when the product of two graphs has an anti-van der Waerden number of four.

\begin{thm}\cite{RSW}\label{theorem:rsw}
    If $G$ and $H$ are connected graphs and $|G|,|H| \ge 2$, then $\aw(G\square H,3) \le 4$.
\end{thm}

The next result, Theorem \ref{PmxPn}, is used when we find an isometric subgraph $P_m\square P_n$ within $G\square H$ to either show that our coloring must have a rainbow within the isometric subgraph or that we know we can color the isometric subgraph in some way and avoid rainbows.

\begin{thm}\cite{RSW}\label{PmxPn}
    For $m,n \geq 2$, \[\aw(P_m \square P_n, 3) = \begin{cases}
3 & \text{if $m = 2$ and $n$ is even, or $m = 3$ and $n$ is odd,} \\
4 & \text{otherwise.}
\end{cases}\]
\end{thm}

Some results in the paper rely on identifying a spine of a tree $T$. A \emph{spine} $S$ of a tree $T$ is defined as a path subgraph of $T$ that realizes $\diam(T)$.  We say $v \in V(S)$ is the \emph{root} of $u \in V(T) \setminus V(S)$ if $v$ is the closest vertex of $S$ to $u$.  
The \emph{branch} $B_v$ of vertex $v$ is the induced subgraph of $T$ containing all vertices whose root is $v$.
Specifically, \[V(B_v) = \{u: \text{$v$ is the root of $u$}\}.\]
Note that $v\in B_v$.

\begin{lem}\label{equidistantspinevertex}
Suppose $T$ is a tree with a spine $S$. If $u \in V(T)\setminus V(S)$ such that $u \in V(B_x)$ and $v \in V(T) \setminus V(B_x)$, then there exists a $u' \in V(S)$ such that $\dist(u,v) = \dist(u',v)$.
\end{lem} 

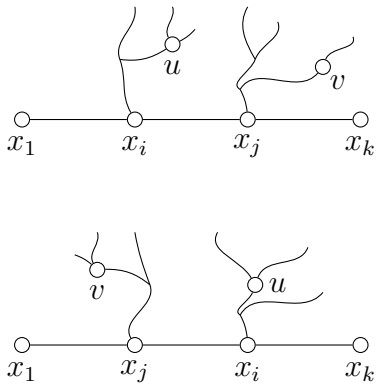
\begin{figure}[h!]
    \centering

    
    \begin{tikzpicture}
        \node[draw,circle,inner sep = 2] (11) at (0,4) {};
        \node[draw,circle,inner sep = 2] (ii) at (1.5,4) {};
        \node[draw,circle,inner sep = 2] (jj) at (3,4) {};
        \node[draw,circle,inner sep = 2] (kk) at (4.5,4) {};
        \node[draw,circle,inner sep = 2] (uu) at (2,5) {};
        \node[draw,circle,inner sep = 2] (vv) at (4,4.7) {};
        
        \draw (11) -- (ii);
        \draw (ii) -- (jj);
        \draw (jj) -- (kk);
        \draw (ii) to[in = -70, out = 120] (1.3,4.8);
        \draw[bend left = 10, bend right = 20] (1.3,4.8) to (uu);
        \draw (1.3,4.8) to[in = -80, out = 110] (1.5,5.5);
        \draw (uu) to[in = -70, out = 120] (2,5.5);
        \draw[bend right = 10] (uu) to (2.3,5.2);
        \draw[bend right = 10] (jj) to (2.9,4.4);
        \draw (2.9,4.4) to[in = -130, out = 60] (vv);
        \draw (vv) to[in = -70, out = 70] (4.4,5.1);
        \draw (2.9,4.4) to[in = -110, out = 150] (3.1,4.8);
        \draw (3.1,4.8) to[in = -70, out = 70] (3.4,5.3);
        \draw (3.1,4.8) to[in = -120, out = 120] (3,5.5);
        
        \node at (0,3.65) {$x_1$};
        \node at (1.5,3.65) {$x_i$};
        \node at (3,3.65) {$x_j$};
        \node at (4.5,3.65) {$x_k$};
        \node at (2,4.7) {$u$};
        \node at (4.2,4.5) {$v$};


        \node[draw,circle,inner sep = 2] (1) at (0,1) {};
        \node[draw,circle,inner sep = 2] (j) at (1.5,1) {};
        \node[draw,circle,inner sep = 2] (i) at (3,1) {};
        \node[draw,circle,inner sep = 2] (k) at (4.5,1) {};
        \node[draw,circle,inner sep = 2] (v) at (1,2) {};
        \node[draw,circle,inner sep = 2] (u) at (3.1,1.8) {};
        
        \draw (1) -- (j);
        \draw (j) -- (i);
        \draw (i) -- (k);
        \draw (j) to[in = -70, out = 120] (1.7,1.8);
        \draw[bend left = 10, bend right = 20] (1.7,1.8) to (v);
        \draw (1.7,1.8) to[in = -80, out = 110] (1.5,2.5);
        \draw (v) to[in = -70, out = 120] (1,2.5);
        \draw[bend right = 10] (v) to (.7,2.2);
        \draw[bend right = 10] (i) to (2.9,1.4);
        \draw (2.9,1.4) to[in = -130, out = 60] (4,1.7);
        \draw (2.9,1.4) to[in = -110, out = 150] (u);
        \draw (u) to[in = -110, out = 70] (3.8,2.3);
        \draw (u) to[in = -120, out = 120] (2.6,2.5);
        
        \node at (0,.65) {$x_1$};
        \node at (1.5,.65) {$x_j$};
        \node at (3,.65) {$x_i$};
        \node at (4.5,.65) {$x_k$};
        \node at (1,1.7) {$v$};
        \node at (3.4,1.8) {$u$};
    \end{tikzpicture}
    \caption{The two cases from the proof of Lemma \ref{equidistantspinevertex}.  The top figure is when $i < j \le k$ and the bottom figure is when $1 \le j < i$.}\label{fig:spine}
\end{figure}

\bpf
Assume $V(S) = \{x_1,\ldots, x_k\}$ with edges $x_ix_{i+1}$ for $1 \leq i \leq k-1$, $u \in V(B_{x_i})$ and $v \in V(B_{x_j})\subset V(T) \setminus V(B_{x_i})$.  
For the sake of contradiction, assume there is no $x_\ell\in V(S)$ with $\dist(x_\ell,v) = \dist(u,v)$.  
First consider the case when $i< j \le k$ (see Figure \ref{fig:spine}).  
This implies that $\dist(x_1,v) < \dist(u,v)$, if not then there would be some $\ell\in \{x_1,x_2,\dots,x_i\}$ with $\dist(x_\ell,v) = \dist(u,v)$.  
However, $\dist(x_1,v) < \dist(u,v)$ implies that $\dist(x_1,x_i) < \dist(u,x_i)$.  
Adding $\dist(x_i,x_k)$ to both sides gives that $\dist(x_1,x_k) < \dist(u,x_k)$ which contradicts that $S$ is a spine.  
A similar contradiction arises if instead it is assumed that $1\le j < i$.  
In this case, $\dist(x_k,v) < \dist(u,v)$.  
This implies that $\dist(x_k,x_1) < \dist(u,x_1)$, again contradicting that $S$ is a spine.  
Thus, there is some $u'\in V(S)$ such that $\dist(u,v) = \dist(u',v)$.
\epf

Many of the conditions in this paper are about peripheral vertices and diameter so some definitions are needed. For a vertex $v$ in a connected graph $G$ the \emph{eccentricity} of $v$, denoted $\epsilon(v)$, is the distance between $v$ and a vertex furthest from $v$ in $G$.  If a vertex $v$ has $\epsilon(v) = \diam(G)$ we call $v$ a \emph{peripheral} vertex. If a graph $G$ contains vertices $u_1,\ldots, u_{n}$ such that $\dist(u_i,u_j) = \diam(G)$ for all distinct $i,j \in \{1,\ldots,n\}$, then we call $G$ \emph{$n$-peripheral}.  A graph is \emph{not $n$-peripheral} if we cannot find $n$ vertices that are pairwise diameter away from each other.
Specifically, we focus on graphs that are $3$-peripheral and graphs that are not $3$-peripheral.

Some of the rainbow-free colorings that are used to get a lower bound on $\aw(G,3)$ depend on using the path between two peripheral vertices that realize $\diam(G)$. Corollary~\ref{d(u,v)+1} provides a distance result for non-peripheral vertices.

\begin{cor}\label{d(u,v)+1}
If $u$ is not a peripheral vertex of $T$ and $v \in V(T)$, then there exists a vertex $w \in V(T)$ such that $\dist(w,v) = \dist(u,v) + 1$. 
\end{cor}

\bpf Let $S$ be a spine of $T$, $V(S) = \{x_1,\dots,x_k\}$ with $x_ix_{i+1}\in E(T)$ for $1\le i \le k-1$.  
If $u \in V(S)$, then, since $u$ is not a peripheral vertex, $u = x_i$ for some $2\le i \le k-1$.  
Thus, either $\dist(x_{i-1},v) = \dist(u,v) + 1$ or $\dist(x_{i+1},v) = \dist(u,v) + 1$.  
If $u\notin V(S)$, Lemma \ref{equidistantspinevertex} indicates there is an $x_i\in V(S)$ such that $\dist(x_i,v) = \dist(u,v)$.  
Since $u$ is not a peripheral vertex neither is $x_i$ so, again, either $\dist(x_{i-1},v) = \dist(u,v) + 1$ or $\dist(x_{i+1},v) = \dist(u,v) + 1$.
\epf

Lemma~\ref{equidistantspinevertex} and Corollary~\ref{d(u,v)+1} are used in Section~\ref{3peripheraltrees} to prove anti-van der Waerden results on the graph products with $3$-peripheral trees.

\section{Trees that are $3$-peripheral}\label{3peripheraltrees}

The idea and motivation for studying $3$-peripheral trees arose from noticing that the path $P_n$ maximizes the diameter and the star $K_{1,n-1}$  minimizes the diameter for all trees of order $n$. Thus, all trees of order $n$ may be considered to be `between' $K_{1,n-1}$ and $P_n$. In working on results associated with stars, we  realized that $\aw(G\square K_{1,n},3)=3$ for all $G$ and $n \ge 3$. The structure of the star allows for this argument,  so we generalized this idea to the class of trees we call $3$-peripheral. 
This led to us showing that $\aw(G\square T,3)=3$ for all connected graphs $G$ with $|G| \ge 2$ and all $3$-peripheral trees $T$ in the main result of this section, Theorem \ref{TxG not $2$-per}.  To begin we present a general result.

\begin{lem}\label{Not$2$-peripheral dist}
    If $T$ is $3$-peripheral, then $\diam(T)$ is even.  Further, for any three vertices that are pairwise distance $\diam(T)$ apart, there is some vertex that is equidistant from all three of them.
\end{lem}

\begin{proof}
    Since $T$ is $3$-peripheral there exist $x,y,z\in V(T)$ with $\dist(x,y) = \dist(y,z) = \dist(x,z) = \diam(T)$.  
    First, assume that the $x-y$ path and $y-z$ path only intersect at $y$.  
    Note that in a tree there is a unique path between any two vertices, thus the $x-z$ path is the path created by traversing the $x-y$ path and then the $y-z$ path.  
    This implies that $\dist(x,z) = 2\diam(T)$, a contradiction. 
    So now assume that the $x-y$ path and the $y-z$ path intersect at a vertex that is not $y$.
    Let $y'$ be the vertex that is closest to $z$ and is on both the $x-y$ path and $y-z$ path. This means
    \[ \dist(x,y) = \dist(x,y') + \dist(y',y) = \diam(T),\]
    \[\dist(z,y) = \dist(z,y') + \dist(y',y) = \diam(T),\] and 
    \[\dist(x,z) = \dist(x,y') + \dist(y',z) = \diam(T).\]

    These equations yield $\dist(y,y') = \dist(x,y') = \dist(z,y')$ so $\diam(T)$ must be even and $y'$ is equidistant from $x,y$ and $z$.
\end{proof}

Lemma \ref{lemma:P2xT} is the base case for Lemma \ref{PnxT not $2$-per}. Lemma \ref{PnxT not $2$-per} is then used to prove the main result of the section, Theorem \ref{TxG not $2$-per}.

\begin{lem}\label{lemma:P2xT}
    If $T$ is a  $3$-peripheral tree, then $\aw(P_2 \square T, 3) = 3$.
\end{lem}

\bpf
Let $c$ be an exact 3-coloring of $P_2 \square T$ and let $T_1$ and $T_2$ denote the two labeled copies of $T$. Assume $c$ avoids rainbow $3$-APs. By Lemma \ref{isometricpathorC3}, there either exists an isometric path or a $C_3$ containing $red$, $blue$, and $green$ in $P_2 \square T$. The latter case produces a rainbow $3$-AP, thus there exists an isometric path in $P_2 \square T$ with $red$, $blue$, and $green$. Let $P$ be the shortest such path, and suppose $P$ is the $v_{j, \alpha} - v_{k, \beta}$ geodesic. Without loss of generality, suppose $c(v_{j, \alpha}) = blue$, $c(v_{k,\beta}) = green$, and $c(v) = red$ for all $v \in V\big(P\setminus \{v_{j,\alpha}. v_{k,\beta}\}\big)$. Lemma \ref{|c(Hi)|<3} implies that $j \neq k$. Without loss of generality, suppose $j = 1, k= 2$.
\par
Assume $\dist(v_{1,\alpha},v_{2,\beta})$ is even. Then there exists a vertex $v$ equidistant from $v_{1,\alpha}$ and $v_{2,\beta}$ on $P$. So, $\{v_{1,\alpha},v,v_{2,\beta}\}$ is a rainbow $3$-AP of common difference $\frac{\dist(v_{1,\alpha},v_{2,\beta})}{2}$.
\par
Thus, $\dist(v_{1,\alpha},v_{2,\beta})$ is odd. Let $S$ denote a spine of $T$. Assume $\dist_T(w_{\alpha},w_{\beta}) < \diam(T)$.
\begin{description}
\item[Case 1.] Suppose $w_{\alpha}, w_{\beta} \in V(S)$. 
\par
Since $\dist_T(w_{\alpha},w_{\beta}) < \diam(T)$, there exists some $w \in V(S)$ such that either \[\dist_T(w,w_{\beta}) = \dist_T(w_{\alpha},w_{\beta}) + 1\] or \[\dist_T(w,w_{\alpha}) = \dist_T(w_{\alpha},w_{\beta}) + 1.\] 
Without loss of generality, suppose $\dist_T(w,w_{\beta}) = \dist_T(w_{\alpha},w_{\beta}) + 1$. 
Let $P$ denote the $w-w_{\beta}$ geodesic in $T$. 
Note that the length of $P$ is odd, and $|c(P_2 \square P)| = 3$. 
Now, Corollary \ref{cor:isosubprod} and Theorem \ref{PmxPn} imply that the isometric subgraph $P_2 \square P$ in $P_2 \square T$ contains a rainbow $3$-AP. 
\item[Case 2.] Suppose $w_{\alpha} \notin V(S)$ or $w_{\beta} \notin V(S)$.
\par
Assume $w_{\alpha},w_{\beta} \notin V(S)$ where $w_{\alpha}$ and $w_{\beta}$ are in the same branch, say $B_x$. 
Without loss of generality, suppose $\dist_T(w_{\beta},x) < \dist_T(w_{\alpha},x)$. 
Then we can find a $w \in B_x \cup \{x\}$ such that $\dist(w_{\alpha},w) = \dist(w_{\alpha},w_{\beta})$. 
Define $P'$ to be the $w_{\alpha}-w$ path in $H$. 
Similar to the first case, Corollary \ref{cor:isosubprod} and Theorem \ref{PmxPn} imply that the isometric subgraph $P_2 \square P'$ in $P_2 \square T$ contains a rainbow $3$-AP. 
\par
If $w_{\alpha} \notin V(S)$ and $w_{\beta}$ isn't in the branch containing $w_{\alpha}$, then, by Lemma \ref{equidistantspinevertex}, there exists some $w_{\alpha'} \in V(S)$ such that $\dist_T(w_{\alpha},w_{\beta}) = \dist_T(w_{\alpha'},w_{\beta})$. 
So, the $3$-AP $\{v_{1,\alpha},v_{2,\beta},v_{1,\alpha'}\}$ implies that $c(v_{1,\alpha'}) = blue$. 
\par
Likewise, if $w_{\beta} \notin V(S)$ and $w_{\alpha}$ isn't in the branch containing $w_{\beta}$, then, by Lemma \ref{equidistantspinevertex}, there exists some $w_{\beta'} \in V(S)$ such that $\dist_T(w_{\alpha},w_{\beta}) = \dist_T(w_{\alpha},w_{\beta'})$. The $3$-AP $\{v_{2,\beta},v_{1,\alpha},v_{2,\beta'}\}$ implies that $c(v_{2,\beta'}) = green$. Case 1 implies that there exists a rainbow $3$-AP in $P_3\square H$.
\end{description}
\par
Thus, $\dist_T(w_\alpha,w_\beta) = \diam(T)$. By Lemma \ref{Not$2$-peripheral dist}, there exists some $w_{\gamma}$ such that $\dist_T(w_\alpha,w_\gamma) = \dist_T(w_\beta,w_\gamma) = \diam(T)$. Define \[S = \big\{\dist(v_{1,\alpha}, v_{1,\beta}), \dist(v_{1,\beta}, v_{2,\beta}), \dist(v_{1,\alpha}, v_{2,\alpha}), \dist(v_{2,\alpha}, v_{2,\beta})\big\}.\]

\begin{stretchpars}Note that $s < \dist(v_{1,\alpha}, v_{2,\beta})$ for all $s \in S$, so the minimality of $\dist(v_{1,\alpha}, v_{2,\beta})$ implies that $c(v_{1,\beta}) = c(v_{2,\alpha}) = red$.
The $3$-APs $\{v_{1,\alpha},v_{2,\beta},v_{1,\gamma}\}$ and\end{stretchpars}\vspace{-.1in} $\{v_{2,\alpha},v_{1,\gamma},v_{2,\beta}\}$ imply that $c(v_{1,\gamma}) = green$, but then $\{v_{1,\alpha}, v_{1,\beta}, v_{1,\gamma}\}$ is a rainbow 3-AP. Thus, $\aw(P_2 \square T, 3) = 3$.
\epf

\begin{lem}\label{PnxT not $2$-per}
    If $n \geq 2$ and $T$ is a $3$-peripheral tree, then $\aw(P_n \square T, 3) = 3$.
\end{lem}

\bpf
Suppose $T$ is a $3$-peripheral tree. 
By Lemma \ref{lemma:P2xT}, $\aw(P_2 \square T, 3) = 3$. 
Let, $\ell \geq 2$ and suppose $\aw(P_{\ell} \square T, 3) = 3$. 
Let $c$ be an exact 3-coloring of $P_{\ell + 1} \square T$, without loss of generality, say $red,blue, green \in c(V(P_{\ell + 1}\square T))$, and assume $c$ avoids rainbow $3$-APs. 
Let $T_i$ denote the $i$th labeled copy of $T$. 
By hypothesis, \[\left|c\left(\bigcup_{i=1}^{\ell}V(T_i)\right)\right| \leq 2 \text{ and } \left|c\left(\bigcup_{i=2}^{\ell+1}V(T_i)\right)\right| \leq 2.\]

Assume $|c\left(\bigcup_{i=2}^{\ell}V(T_i)\right)| = 2$, say $red,blue \in c\left(\bigcup_{i=2}^{\ell}V(T_i)\right)$. 
Lemma \ref{|c(V(Gi U Gj))|<3} implies that $c(V(T_1)) \cup c(V(T_{\ell+1})) \subseteq \{red,blue\}$, contradicting that $green \in c(V(P_{\ell+1} \square T))$. 
So, without loss of generality, suppose \[c\left(\bigcup_{i=2}^{\ell}V(T_i)\right) = \{red\}.\]

If $c(V(T_1)) = \{red\}$, then $blue, green \in c(V(T_{\ell+1}))$, contradicting Lemma \ref{|c(V(Gi U Gj))|<3}. 
If $blue, green \in c(V(T_1))$, then we contradict Lemma \ref{|c(V(Gi U Gj))|<3} again. 
Thus, without loss of generality, $blue \in c(V(T_1))$ and $green \in c(V(T_{\ell+1}))$. 
Say $c(v_{1,j}) = blue$ and $c(v_{\ell+1,k}) = green$, and without loss of generality, suppose $j \leq k$. 
\par 
Assume $\dist(w_j,w_k) < \diam(T)$ and let $S$ denote a spine of $T$. First, suppose that, without loss of generality, $w_j$ is not a peripheral vertex. Corollary \ref{d(u,v)+1} implies there exists a $w_i \in V(T)$ such that $\dist(w_i,w_k) = \dist(w_j,w_k) + 1$. Otherwise, if both $w_j$ and $w_k$ are peripheral vertices, then there exists some vertex $w_\ell$ such that $\dist(w_k,w_\ell) = \diam(T)$. Since $\dist(w_j,w_k) < \diam(T)$, there exists some $w_i$ on the $w_k-w_\ell$ path such that $\dist(w_i,w_k) = \dist(w_j,w_k) + 1$. Now, Lemma \ref{|c(V(Gi U Gj))|<3} and the 3-AP $\{v_{1,j},v_{\ell+1,k},v_{2,i}\}$ imply that $c(v_{2,i}) = blue$. However, this contradicts that $c(T_2) = \{red\}$. Thus, $\dist(w_j,w_k) = \diam(T)$.
\par
By Lemma \ref{Not$2$-peripheral dist}, there exists some $w_x$ such that $\dist_T(w_j,w_x) = \dist_T(w_k,w_x) = \diam(T)$. Define \[S = \big\{\dist(v_{1,j}, v_{1,k}), \dist(v_{1,k}, v_{\ell+1,k}), \dist(v_{1,j}, v_{\ell+1,j}), \dist(v_{\ell+1,j}, v_{\ell+1,k})\big\}.\] Note that $s < \dist(v_{1,j}, v_{\ell+1,k})$ for all $s \in S$. The minimality of $\dist(v_{1,j}, v_{\ell+1,k})$ implies that $c(v_{1,k}) = c(v_{\ell+1,j}) = red$.
The $3$-APs $\{v_{1,j},v_{\ell+1,k},v_{1,x}\}$ and $\{v_{\ell+1,j},v_{1,x},v_{\ell+1,k}\}$ imply that $c(v_{1,x}) = green$, but then $\{v_{1,j}, v_{1,k}, v_{1,x}\}$ is a rainbow 3-AP. Thus, $\aw(P_{\ell+1} \square T, 3) = 3$. By induction, $\aw(P_n \square T, 3) = 3$ for all $n \geq 2$.
\epf

\begin{thm}\label{TxG not $2$-per}
    If $T$ is a $3$-peripheral tree and $G$ is connected with $2\le |G|$, then \[\aw(T \square G, 3) = 3.\]
\end{thm}

\bpf
Let $c$ be an exact 3-coloring of $T \square G$ where $|G| = \{w_1,\ldots,w_n\}$ and let $T_i$ denote the $i$th labeled copy of $T$. 
Assume $c$ avoids rainbow $3$-APs. 
Let $\alpha = \max_{1\le i \le n}\{|c(T_i)|\}$ and choose $k$ such that $|c(T_k)| = \alpha$. 
If $\alpha \leq 2$, there is a path of at most length two in $T_k$ containing $\alpha$ colors. 
If $\alpha = 3$, Lemma \ref{isometricpathorC3} implies that $T_k$ contains an isometric $C_3$ or path containing all three colors. If there exists an isometric $C_3$ with all three colors, then there exists an immediate rainbow 3-AP in $T_k$. 
Thus, there must exist an isometric path in $T_k$ containing $\alpha$ colors. 
Let $P$ be a shortest such path.
\begin{description}
\item[Case 1.] Suppose $P$ contains three colors.
\par
First suppose the length of $P$ is even, and, without loss of generality, suppose the two leaves of $P$ are colored $red$ and $blue$ and the rest of the vertices are colored $green$. Then there exists a $green$ vertex equidistant from the $red$ and $blue$ ones, creating a rainbow $3$-AP.
\par
Second suppose the length of $P$ is odd, and $T_k$ contains three colors. 
Note that $w_i \in V(G)$ has a neighbor, say $w_j$. 
Then let $P_2$ denote the $w_i -w_j$ path in $G$, and let $P^*$ denote the path corresponding to $P$ in $T$. 
Then the isometric subgraph $P_2 \square P^*$ of $G\square H$ contains a rainbow $3$-AP by Corollary \ref{cor:isosubprod} and Theorem \ref{PmxPn}.
\item[Case 2.] Suppose $P$ contains one color.
\par
This implies that each $T_k$ is monochromatic by the definition of $P$. 
Since $T \square G$ has three colors, there exists either an isometric $C_3$ or an isometric shortest path $\rho$ in a
copy of $G$ that has at least three colors by Lemma \ref{isometricpathorC3}. If there is an isometric $C_3$, then there is an immediate rainbow $3$-AP. 
In the other case, this is just Case 1 with the roles of $T$ and $G$ reversed.
\item[Case 3.] Suppose $P$ contains two colors.
\par
Then there must exist some $T_k$ with two colors, say $T_k = \{red,blue\}$. 
Lemma \ref{c(H_i)/c(H_j) < 2} implies that $green$ must appear either with $red$ or $blue$, without loss of generality,
suppose $green$ appears with $red$.
\par
Define $c'$ to be the auxiliary coloring of $G$ where \[c'(w_{i}) = \begin{cases}
red & \text{if $c(V(T_i)) = \{red\}$,} \\
\C & \text{if $\C \in c(V(T_i))$.}
\end{cases}\]
By Lemma \ref{isometricpathorC3}, there either exists an isometric path containing $red$, $blue$, and $green$ or a $C_3$ containing $red$, $blue$, and $green$. 
The latter case implies an immediate rainbow $3$-AP, thus there exists an isometric path $P'$ in $G$ with $red$, $blue$, and $green$. 
Then by Corollary \ref{cor:isosubprod} and Lemma \ref{PnxT not $2$-per}, there exists a rainbow $3$-AP in the isometric subgraph $T \square P'$ of $T \square G$, contradicting that $c$ avoids rainbow $3$-APs.
\end{description}
In every case, there exists a rainbow $3$-AP in $T \square G$, thus $\aw(T \square G, 3) = 3$.
\epf

\section{Trees that are not $3$-peripheral}\label{sec:not3peripheral}

Section \ref{3peripheraltrees} culminates with a result about graph products between $3$-peripheral trees and arbitrary nontrivial connected graphs. 
This is due to the strong structure that a graph being $3$-peripheral provides. 
The structure of trees which are not $3$-peripheral is much less impactful, and this is discussed in more detail in section \ref{sec:futurework}. 
That is not to say such trees have no harnessable structure as shown by Lemma \ref{lemma:treediam}. 
This section concludes with Theorem \ref{theorem:treeprododd} which provides a class of products of trees which are not $3$-peripheral whose anti-van der Waerden number is four.

\begin{lem}\label{claim: cross product 2-peripheral}
If $G_1$ and $G_2$ are not $3$-peripheral, then $G_1\square G_2$ is not $3$-peripheral.
\end{lem}

\bpf
    We will prove the contrapositive, so assume that $G_1\square G_2$ is $3$-peripheral.
    Then there exist $x,y,z\in V(G_1\square G_2)$ with \[\dist_{G_1\square G_2}(x,y) = \dist_{G_1\square G_2}(y,z) = \dist_{G_1\square G_2}(x,z) = \diam(G_1) + \diam(G_2)\]
    by definition and Observation \ref{obs:diamGxH}.
    Assume that $x,y$ and $z$ correspond to $x_1,y_1,z_1\in V(G_1)$ and $x_2,y_2,z_2\in V(G_2)$, respectively.  
    Now, \[\dist_{G_1\square G_2}(x,y) = \dist_{G_1}(x_1,y_1) + \dist_{G_2}(x_2,y_2), \dist_{G_1\square G_2}(x,z) = \dist_{G_1}(x_1,z_1) + \dist_{G_2}(x_2,z_2) \] and \[\dist_{G_1\square T_2}(y,z) = \dist_{G_1}(y_1,z_1) + \dist_{G_2}(y_2,z_2).\]  
    Thus, $x_1,y_1, z_1\in V(G_1)$ are pairwise distance $\diam(G_1)$ and $x_2,y_2, z_2\in V(G_2)$ are pairwise distance $\diam(G_2)$ away from each other so $G_1$ and $G_2$ are $3$-peripheral.
\epf

Lemma~\ref{lemma:treediam} is vital in determining when a product involving a tree has anti-van der Waerden number four.

\begin{lem}\label{lemma:treediam}
    Suppose $T$ is a tree that is not $3$-peripheral and $u_i,u_j\in V(T)$ realize the diameter of $T$.
    If there exist $u_x,u_y \in V(T)$ such that $\dist(u_x,u_j) = \diam(T)$ and $\dist(u_i,u_y) = \diam(T)$, then $\dist(u_x,u_y) = \diam(T)$. 
\end{lem}

\begin{proof}
    If $u_i  = u_x$ or $u_j = u_y$, then the result is immediate. 
    Define $P_x$ to be the $u_x-u_j$ path, $P_y$ to be the $u_i-u_y$ path, and $P$ to be the $u_i-u_j$ path. 
    Let $u_\ell$ be the closest vertex in $V(P) \cap V(P_x)$ to $u_i$, and $u_k$ to be closest vertex in $V(P) \cap V(P_y)$ to $u_j$. 
    Note that \[\dist(u_x,u_\ell) + \dist(u_\ell,u_j) = \dist(u_x,u_j) = \dist(u_i,u_j) = \dist(u_i,u_\ell) + \dist(u_\ell,u_j).\] 
    Subtracting $\dist(u_{\ell},u_j)$ shows that $\dist(u_x,u_\ell) = \dist(u_i,u_\ell)$. 
    Similarly, $\dist(u_k,u_y) = \dist(u_k,u_j)$. 

    Assume $d(u_i,u_k) < d(u_i,u_\ell)$. 
    Then $P$ is the concatenation of the $u_i - u_k$, $u_k - u_\ell$, and $u_\ell - u_j$ paths. 
    Notice that $\diam(T) = d(u_i,u_y) = d(u_i,u_k) + d(u_k,u_y)$ implies that one of $d(u_i,u_k)$ and $d(u_k,u_y)$ is at least $\diam(T)/2$. 
    Similarly, at least one of $d(u_x,u_\ell)$ and $d(u_\ell,u_j)$ is at least $\diam(T)/2$. 
    Concatenating these two paths with the nontrivial $u_k-u_\ell$ path yields a path whose length is greater than $\diam(T)$, a contradiction. 

    Next, assume $u_k = u_\ell$. 
    If $d(u_i,u_k) \neq \diam(T)/2 \neq d(u_k,u_j)$, then one of the paths $u_i - u_k - u_x$ or $u_y - u_k - u_j$ has length greater than $\diam(T)$, a contradiction. 
    Thus, $d(u_i,u_k) = d(u_k,u_j) = \diam(T)/2$. 
    Since \[\diam(T) = d(u_i,u_y) = d(u_i,u_k) + d(u_k,u_y) = \frac{\diam(T)}{2} + d(u_k,u_y),\] it follows that $d(u_k,u_y) = \diam(T)/2$. 
    Since $u_k = u_\ell$ lies on the $u_y-u_j$ path, we have $d(u_y,u_j) = d(u_y,u_k) + d(u_k,u_j) = \diam(T)$. 
    Now $u_i,u_j,u_y$ pairwise realize the diameter, contradicting that $T$ is not $3$-peripheral.
    
    Finally, suppose $d(u_i,u_\ell) < d(u_i,u_k)$ so that $P$ is the concatenation of the $u_i - u_\ell$, $u_\ell - u_k$, and $u_k - u_j$ paths. 
    Observe, 
    
    \begin{equation*}
        \begin{split}
            \dist(u_x,u_y) & = \dist(u_x,u_\ell) + \dist(u_\ell,u_k) +\dist(u_k,u_y) \\
            & = \dist(u_i,u_\ell) + \dist(u_\ell,u_k) +\dist(u_k,u_j) \\
            & = \dist(u_i,u_j) \\
            & = \diam(T).
        \end{split} 
    \end{equation*}
\end{proof}

\begin{thm}\label{theorem:treeprododd}
If $T$ and $T'$ are trees which are not $3$-peripheral with $|T|,|T'| \geq 2$ and $\diam(T \square T')$ is odd, then $\aw(T \square T', 3) = 4$. 
\end{thm}

\bpf
Suppose $V(T) = \{u_1,\ldots, u_x\}$ and $V(T') = \{w_1,\ldots, w_y\}$. Let $S$ and $S'$ be the spines defined by the $u_i-u_j$ and $w_h-w_k$ geodesics, respectively.  Define $c: V(T \square T') \to \{red,blue,green\}$ by \[c(v_{\alpha,\beta}) = 
\begin{cases}
    red & \text{if  
$\dist(v_{\alpha,\beta},v_{j,k}) = \diam(T \square T')$,} \\
    blue & \text{if $\dist(v_{\alpha,\beta},v_{i,h}) = \diam(T \square T')$,} \\
    green & \text{otherwise.}
\end{cases}\]
By Lemma \ref{claim: cross product 2-peripheral}, $T \square T'$ is not $3$-peripheral and so $c$ is a well-defined coloring. Note that if a rainbow $3$-AP exists, it must contain a $red$ and a $blue$ vertex. Let $v_{r_1,r_2}$ and $v_{b_1,b_2}$ denote arbitrary $red$ and $blue$ vertices, respectively. Proposition \ref{proposition:dist} and Observation \ref{obs:diamGxH} imply that $\dist(u_{r_1},u_j)=\diam(T) = \dist(u_{b_1},u_i)$ and $\dist(w_{r_2},w_k)=\diam(T') = \dist(w_{b_2},w_h)$. Now, Lemma \ref{lemma:treediam} implies that $\dist(u_{r_1},u_{b_1}) = \diam(T)$ and $\dist(w_{r_2},w_{b_2})=\diam(T')$, so \[\dist(v_{r_1,r_2},v_{b_1,b_2}) = \diam(T \square T').\]

\begin{stretchpars}Now, for the sake of contradiction, assume there exists some $v_{\alpha,\beta} \in V(T\square T')$ such that $\{v_{b_1,b_2},v_{r_1,r_2},v_{\alpha,\beta}\}$ or $\{v_{r_1,r_2},v_{b_1,b_2},v_{\alpha,\beta}\}$ is a $3$-AP. 
If\end{stretchpars}\vspace{-.1in} $\{v_{b_1,b_2},v_{r_1,r_2},v_{\alpha,\beta}\}$ is a $3$-AP, then $\dist(v_{r_1,r_2},v_{\alpha,\beta}) = \diam(T\square T')$. 
Thus, Proposition \ref{proposition:dist} and Observation \ref{obs:diamGxH} imply that 

\[\dist(u_{r_1},u_\alpha)=\diam(T) \text{ and }\dist(w_{r_2},w_\beta)=\diam(T').\]

The first equality, along with $\dist(u_i,u_{b_1}) = \diam(T)$ and Lemma \ref{lemma:treediam}, imply that $\dist(u_i,u_\alpha)=\diam(T)$. 
Similarly, the second equality, along with $\dist(w_h,u_{b_2}) = \diam(T')$ and Lemma \ref{lemma:treediam}, imply that $\dist(w_h,w_\beta)=\diam(T')$. 
Thus, $c(v_{\alpha,\beta})=blue$, and $\{v_{b_1,b_2},v_{r_1,r_2},v_{\alpha,\beta}\}$ is not a rainbow $3$-AP.  
A similar argument shows that $\{v_{r_1,r_2},v_{b_1,b_2},v_{\alpha,\beta}\}$ is not a rainbow $3$-AP.

Thus, any rainbow $3$-AP must be of the form $\{v_{r_1,r_2},v_{\alpha,\beta},v_{b_1,b_2}\}$ for some $v_{\alpha,\beta} \in V(T \square T')$ meaning $\dist(v_{r_1,r_2},v_{\alpha,\beta}) = \dist(v_{\alpha,\beta},v_{b_1,b_2})$.

 Since $T\square T'$ is bipartite, we can partition $V(T\square T')$ into partite sets $A$ and $B$. 
 Because $\dist(v_{r_1,r_2},v_{b_1,b_2})$ is odd, the vertices $v_{r_1,r_2}$ and $v_{b_1,b_2}$ are in different parts, say $v_{r_1,r_2}\in A$ and $v_{b_1,b_2}\in B$. 
 If $\dist(v_{r_1,r_2},v_{\alpha,\beta})$ is even, then $v_{\alpha,\beta}\in A$ and therefore $\dist(v_{b_1,b_2},v_{\alpha,\beta})$ is odd.
 Similarly, if $\dist(v_{r_1,r_2},v_{\alpha,\beta})$ is odd, then $\dist(v_{b_1,b_2},v_{\alpha,\beta})$ is even. 
So, no such $v_{\alpha,\beta}$ exists, and $c$ contains no rainbow $3$-APs. 
Therefore, $\aw(T \square T', 3) \geq 4$, and Theorem \ref{theorem:rsw} implies that $\aw(T \square T', 3) = 4$.
\epf

\section{Future Work}\label{sec:futurework}

In Section \ref{3peripheraltrees}, $\aw(T\square G,3)$ was determined when $T$ is $3$-peripheral.  In Section \ref{sec:not3peripheral}, $\aw(T\square T',3)$ was determined when $T$ and $T'$ are not $3$-peripheral and $\diam(T\square T')$ is odd. 
The remaining case to consider is when $T$ and $T'$ are not $3$-peripheral and $\diam(T\square T')$ is even.

To illustrate why this case is more difficult, recall that Theorem \ref{PmxPn} implies that if $\diam(P_m\square P_n)$ is even, then $\aw(P_m\square P_n,3)=3$ only if $m\leq3$ or $n\leq3$. When $\diam(P_m\square P_n)$ is odd the coloring which colors $v_{1,1}$ red, $v_{m,n}$ blue, and the remaining vertices green avoids rainbow $3$-APs. 
This was generalized in Section $3$ to the coloring found in the proof of Theorem \ref{theorem:treeprododd}.
Alternatively, when $\diam(P_m\square P_n)$ is even with $m,n\geq4$ the coloring which colors $v_{1,1}$ red, $v_{m-1,n}$ and $v_{m,n-1}$ blue, and the remaining vertices green avoids rainbow $3$-APs. 
When attempting to generalize this coloring, problems arise. 
Following the set-up in the proof of Theorem \ref{theorem:treeprododd}, in particular $v_{j,k}$ and $v_{i,h}$ are diameter vertices in $T\square T'$, the generalized coloring $c:V(T\square T') \to \{red,blue,green\}$ can be defined by 

\[c(v_{\alpha,\beta}) = 
\begin{cases}
    red & \text{if $\dist(v_{\alpha,\beta},v_{j,k}) = \diam(T \square T')$,} \\
    blue & \text{if $\dist(v_{\alpha,\beta},v_{i,h}) = \diam(T \square T')-1$,} \\
    green & \text{otherwise.}
\end{cases}
\]

\begin{figure}[h!]
    \centering
    \begin{tikzpicture}[scale = .7]
        \node[draw, very thick, circle, color=cbgreen] (11) at (0,4.5) {};
        \node[draw,circle,fill=cbred,thick] (12) at (0,3) {};
        \node[draw,circle,thick] (13) at (0,1.5) {};
        \node[draw,circle,thick] (14) at (-.5,0) {};
        \node[draw,circle,thick] (15) at (.5,0) {};
        
        \node[draw,circle,fill=cbred,thick] (21) at (2,4.5) {};
        \node[draw,circle,thick] (22) at (2,3) {};
        \node[draw,circle,thick] (23) at (2,1.5) {};
        \node[draw,circle,thick] (24) at (1.5,0) {};
        \node[draw,circle,thick] (25) at (2.5,0) {};
        
        \node[draw,circle,thick] (31) at (4,4.5) {};
        \node[draw,circle,thick] (32) at (4,3) {};
        \node[draw,circle,thick] (33) at (4,1.5) {};
        \node[draw, very thick, circle, color=cbcyan] (34) at (3.5,0) {};
        \node[draw, very thick, circle, color=cbcyan] (35) at (4.5,0) {};
        
        \node[draw,circle,thick] (41) at (6,4.5) {};
        \node[draw,circle,thick] (42) at (6,3) {};
        \node[draw, very thick, circle, color=cbcyan] (43) at (6,1.5) {};
        \node[draw,circle,fill=cbyellow,thick] (44) at (5.5,0) {};
        \node[draw,circle,fill=cbyellow,thick] (45) at (6.5,0) {};

        \node[draw,circle,fill=cbred,thick] (67) at (-0.5,0) {};
        
        \draw[thick]  (11) to node [auto] {} (12);
        \draw[thick]  (12) to node [auto] {} (13);
        \draw[thick]  (13) to node [auto] {} (14);
        \draw[thick]  (13) to node [auto] {} (15);
        
        \draw[thick]  (21) to node [auto] {} (22);
        \draw[thick]  (22) to node [auto] {} (23);
        \draw[thick]  (23) to node [auto] {} (24);
        \draw[thick]  (23) to node [auto] {} (25);
        
        \draw[thick]  (31) to node [auto] {} (32);
        \draw[thick]  (32) to node [auto] {} (33);
        \draw[thick]  (33) to node [auto] {} (34);
        \draw[thick]  (33) to node [auto] {} (35);
        
        \draw[thick]  (41) to node [auto] {} (42);
        \draw[thick]  (42) to node [auto] {} (43);
        \draw[thick]  (43) to node [auto] {} (44);
        \draw[thick]  (43) to node [auto] {} (45);
        
        \draw[thick]  (11) to node [auto] {} (21);
        \draw[thick]  (12) to node [auto] {} (22);
        \draw[thick]  (13) to node [auto] {} (23);
        \draw[thick, bend right = 30]  (14) to node [auto] {} (24);
        \draw[thick, bend right = 30]  (15) to node [auto] {} (25);
        
        \draw[thick]  (31) to node [auto] {} (21);
        \draw[thick]  (32) to node [auto] {} (22);
        \draw[thick]  (33) to node [auto] {} (23);
        \draw[thick, bend left = 30]  (34) to node [auto] {} (24);
        \draw[thick, bend left = 30]  (35) to node [auto] {} (25);
        
        \draw[thick]  (31) to node [auto] {} (41);
        \draw[thick]  (32) to node [auto] {} (42);
        \draw[thick]  (33) to node [auto] {} (43);
        \draw[thick, bend right = 30]  (34) to node [auto] {} (44);
        \draw[thick, bend right = 30]  (35) to node [auto] {} (45);

        \node[draw, thick, circle] (11) at (9.5,4.5) {};
        \node[draw,circle,thick] (12) at (9,3) {};
        \node[draw,circle,thick] (13) at (9,1.5) {};
        \node[draw,circle,thick] (14) at (8.5,0) {};
        \node[draw,circle,thick] (15) at (9.5,0) {};
        \node[thick, draw, circle] (16) at (8.5,4.5) {};
        
        \node[draw,circle,thick] (21) at (11.5,4.5) {};
        \node[draw,circle,thick] (22) at (11,3) {};
        \node[draw,circle,thick] (23) at (11,1.5) {};
        \node[draw,circle,thick] (24) at (10.5,0) {};
        \node[draw,circle,thick] (25) at (11.5,0) {};
        \node[thick, draw, circle] (26) at (10.5,4.5) {};
        
        \node[draw,circle,thick] (31) at (13.5,4.5) {};
        \node[draw,circle,thick] (32) at (13,3) {};
        \node[draw,circle,thick] (33) at (13,1.5) {};
        \node[draw, thick, circle] (34) at (12.5,0) {};
        \node[draw, thick, circle] (35) at (13.5,0) {};
        \node[draw,circle,thick] (36) at (12.5,4.5) {};
        
        \node[draw,circle,thick] (41) at (15.5,4.5) {};
        \node[draw,circle,thick] (42) at (15,3) {};
        \node[draw, thick, circle] (43) at (15,1.5) {};
        \node[draw,circle,thick] (44) at (14.5,0) {};
        \node[draw,circle,thick] (45) at (15.5,0) {};
        \node[draw,circle,thick] (46) at (14.5,4.5) {};
        
        \draw[thick]  (11) to node [auto] {} (12);
        \draw[thick]  (12) to node [auto] {} (13);
        \draw[thick]  (13) to node [auto] {} (14);
        \draw[thick]  (13) to node [auto] {} (15);
        \draw[thick]  (12) to node [auto] {} (16);
        
        \draw[thick]  (21) to node [auto] {} (22);
        \draw[thick]  (22) to node [auto] {} (23);
        \draw[thick]  (23) to node [auto] {} (24);
        \draw[thick]  (23) to node [auto] {} (25);
        \draw[thick]  (22) to node [auto] {} (26);
        
        \draw[thick]  (31) to node [auto] {} (32);
        \draw[thick]  (32) to node [auto] {} (33);
        \draw[thick]  (33) to node [auto] {} (34);
        \draw[thick]  (33) to node [auto] {} (35);
        \draw[thick]  (32) to node [auto] {} (36);
        
        \draw[thick]  (41) to node [auto] {} (42);
        \draw[thick]  (42) to node [auto] {} (43);
        \draw[thick]  (43) to node [auto] {} (44);
        \draw[thick]  (43) to node [auto] {} (45);
        \draw[thick]  (42) to node [auto] {} (46);
        
        \draw[thick, bend left = 30]  (11) to node [auto] {} (21);
        \draw[thick]  (12) to node [auto] {} (22);
        \draw[thick]  (13) to node [auto] {} (23);
        \draw[thick, bend right = 30]  (14) to node [auto] {} (24);
        \draw[thick, bend right = 30]  (15) to node [auto] {} (25);
        \draw[thick, bend left = 30]  (16) to node [auto] {} (26);
        
        \draw[thick, bend right = 30]  (31) to node [auto] {} (21);
        \draw[thick]  (32) to node [auto] {} (22);
        \draw[thick]  (33) to node [auto] {} (23);
        \draw[thick, bend left = 30]  (34) to node [auto] {} (24);
        \draw[thick, bend left = 30]  (35) to node [auto] {} (25);
        \draw[thick, bend right = 30]  (36) to node [auto] {} (26);
        
        \draw[thick, bend left = 30]  (31) to node [auto] {} (41);
        \draw[thick]  (32) to node [auto] {} (42);
        \draw[thick]  (33) to node [auto] {} (43);
        \draw[thick, bend right = 30]  (34) to node [auto] {} (44);
        \draw[thick, bend right = 30]  (35) to node [auto] {} (45);
        \draw[thick, bend left = 30]  (36) to node [auto] {} (46);

         \node[draw,circle,fill=cbred,thick] (12) at (9,3) {};

         \node[draw,circle,fill=cbred,thick] (21) at (11.5,4.5) {};
         \node[draw,circle,fill=cbred,thick] (57) at (8.5,0) {};

         \node[thick, draw, circle, fill=cbred] (26) at (10.5,4.5) {};

         \node[draw,circle,fill=cbyellow,thick] (44) at (14.5,0) {};
         \node[draw,circle,fill=cbyellow,thick] (45) at (15.5,0) {};

    \end{tikzpicture}

    \caption{When we apply the coloring to the product on the left we get two diffrent colorings depending on our choice of $v_{j,k}$ and $v_{i,h}$.  Due to symmetry, there is only one coloring for the graph on the right.}

    \label{fig:cex}
\end{figure}
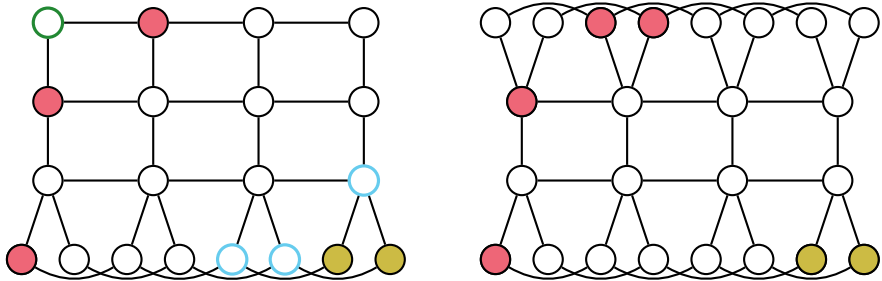

In Figure \ref{fig:cex}, we consider two different products of trees.  
Note that each tree in each product is not $3$-peripheral and the product has even diameter.  The figure on the left has two different colorings based on our choices of $v_{i,h}$ and $v_{j,k}$. 
One coloring uses $\{ {\color{cbred} red}, {\color{cbyellow} yellow}, \mathbbm{white} \}$ and is not rainbow-free and the other coloring uses $\{ {\color{cbgreen} green}, {\color{cbcyan} cyan}, black \}$ and is rainbow-free. 
However, the product on the right has more symmetry so our choice of $v_{i,h}$ and $v_{j,k}$ does not matter and this coloring does have a rainbow.

These examples suggest two things. 
First, $3$-peripheral or not $3$-peripheral alone is not going to fully characterize which tree products have anti-van der Waerden number three and which have anti-van der Waerden number four. 
Second, the product from Figure \ref{fig:cex} on the left was rainbow-free with one coloring and had a rainbow with another coloring. 
This means, there might be a `good choice' required when trying to find a rainbow using a generalized coloring when the diameter of the product is even.

In studying peripheral trees, Theorem \ref{TxG not $2$-per} found that the $\aw(G\square T, 3) = 3$ when $|G| \ge 2$ and $T$ is a $3$-peripheral tree. This led to the following conjecture.

\begin{conj}
    If $T$ is a $k$-peripheral tree and $G$ is a nontrivial connected graph, then $\aw(T\square G,k)=k$.
\end{conj}

\section*{Acknowledgements} 
Thank you to the University of Wisconsin-La Crosse (UWL) Deans Distinguished Fellows program that supported the second and fifth authors.  Thanks also to UWL Undergraduate Research and Creativity grants that supported the second author. Thanks also to the UWL Department of Mathematics and Statistics Bange/Wine Undergraduate Research Endowment that supported the second author.  Finally, thanks to the Simons Laufer Mathematical Sciences Institute (SLMath and formerly MSRI) for supporting the first, third and fourth authors.

\end{document}